\documentstyle[11pt]{article}
\def\Cox{\hfill \Box}

\begin{document}

\title{On near-critical and dynamical percolation in the tree case}
\author{Olle H\"{a}ggstr\"{o}m\thanks{Research supported by a grant from the
Swedish Natural Science Research Council.}
\and Robin Pemantle\thanks{Research supported in part by a Presidential
Faculty Fellowship.} \and \\
Chalmers University of Technology and
University of Wisconsin}
\maketitle

\begin{abstract}
Consider independent bond percolation with retention probability $p$ on a
spherically symmetric tree $\Gamma$. Write $\theta_\Gamma(p)$ for the
probability that the root is in an infinite open cluster, and define the
critical value $p_c=\inf\{p:\, \theta_\Gamma(p)>0\}$. If 
$\theta_\Gamma(p_c)=0$,
then the root may still percolate in the corresponding dynamical percolation
process at the critical value $p_c$, as demonstrated 
recently by H\"aggstr\"om, Peres and Steif. Here we relate this phenomenon
to the near-critical behaviour of $\theta_\Gamma(p)$ by showing
that the root percolates in the dynamical percolation process if and only if 
$\int_{p_c}^1 (\theta_\Gamma(p))^{-1}dp<\infty$. 
The ``only if'' direction extends to general trees, 
whereas the ``if'' direction fails in this generality.  
\end{abstract}

\section{Introduction}

The main setup of this paper is that of independent bond percolation on an 
infinite locally finite connected graph $G=(V,E)$ with a distinguished vertex 
$\rho\in V$. Each edge $e\in E$ is independently assigned value $1$ (open, on) 
with probability $p\in [0,1]$ and $0$ (closed, off)
with probability $1-p$; write
${\bf P}_{G,p}$ for the induced probability measure on $\{0,1\}^E$. Also write 
$\rho\leftrightarrow\infty$ for the event that $\rho$ is in an infinite
open cluster, and define the percolation function
$\theta_G(p)={\bf P}_{G,p}[\rho\leftrightarrow\infty]$. 
Furthermore define the critical value $p_c=p_c(G)=\inf\{p:\, \theta_G(p)>0\}$.
Let ${\cal C}$ denote the event that an infininte open cluster exists
somewhere in $G$. It is well known that $P_{G,p}[{\cal C}]=0$ for $p<p_c$
and $P_{G,p}[{\cal C}]=1$ for $p>p_c$. At $p=p_c$, this 
probability can be either $0$ or $1$, depending on whether 
$\theta_G(p_c)$ is $0$ or strictly positive. This model has been studied
extensively for several decades, mainly in the case where $G$ is the
cubic lattice in $d\geq 2$ dimensions; see e.g. Grimmett \cite{G}. 

We shall also consider a recent extension of the above setup, known as
dynamical percolation. This is a time-dynamical variant which
was introduced by H\"aggstr\"om, Peres
and Steif \cite{HPS}; see also H\"aggstr\"om \cite{H} for a short survey. 
In this model, all edges $e\in E$ turn on and off
according to independent stationary continuous time two-state Markov chains,
turning on at rate $p$ and off at rate $1-p$ (that is, an edge in state $0$
(resp.\ $1$) waits for an amount of time whose distribution is exponential
with mean $\frac{1}{p}$ (resp.\ $\frac{1}{1-p}$) before flipping to state $1$
(resp.\ $0$)). We let $\Psi_{G,p}$ be a probability measure supporting such a 
$\{0,1\}^E$-valued Markov process. The stationary distribution for each
edge puts probabilities $1-p$ and $p$ on states $0$
and $1$, respectively. Hence, by the stationarity assumption, what we see 
at any fixed time $t$ is 
ordinary bond percolation with parameter $p$. Writing ${\cal C}_t$ for the
event that an infinite open cluster exists at time $t$, we thus have (by
Fubini's Theorem) that
\begin{equation} \label{eq:Leb}
\begin{array}{lll}
\Psi_{G,p}[\neg{\cal C}_t \mbox { occurs for a.e.\ $t$}] 
= 1 & \mbox{if} & P_p[{\cal C}]=0
\end{array}
\end{equation}
and similarly
\[
\begin{array}{lll}
\Psi_{G,p}[{\cal C}_t \mbox { occurs for a.e.\ $t$}] 
= 1 & \mbox{if} & P_p[{\cal C}]=1. 
\end{array}
\]
[Here ``a.e.'' is short for ``(Lebesgue-)almost every''. In other words,
the left hand expression in (\ref{eq:Leb}) says that with
$\Psi_{G,p}$-probability $1$, the set of times $t$ for which
$\neg{\cal C}_t$ does {\em not} occur, has Lebesgue measure $0$.]
An obvious question to ask is whether the quantifier ``a.e.'' can be 
strengthened to ``every'' in the above statements. For $p\neq p_c$, the
answer is yes (see \cite{HPS}), but (perhaps surprisingly) the answer 
is no for $p=p_c$ and certain choices of $G$. One of the main problems
in dynamical percolation is the classification of graphs according to whether
exceptional times exist at criticality. 

Here we are interested in graphs for which $\theta_G(p_c)=0$ and in various
notions of how ``close'' such a graph is to percolating at criticality.
We can think of at least two reasonable ways of making ``close'' more
precise.
\begin{enumerate}
\item
The near-critical behaviour of $\theta_G(p)$. How rapidly does $\theta_G(p)$
take off from criticality? For instance, one can ask whether the derivative
$\lim_{p\searrow p_c}\frac{\theta_G(p)}{p-p_c}$ is finite or infinite. 
(Questions of this kind have been studied extensively in terms of so-called
critical exponents, see e.g.\ Chapter 7 of Grimmett \cite{G}.) 
\item
Dynamical percolation. Will $\Psi_{G,p_c}$ assign positive probability to
the existence of exceptional times at which an infinite open cluster exists?
\end{enumerate}
Intuitively, one would expect that if $G$ is sufficiently close to 
percolating at criticality in one of these senses, then it should be close
also in the other. The first result of this kind was obtained
in \cite{HPS}, and says that if 
\begin{equation} \label{derivative}
\limsup_{p\searrow p_c}\frac{\theta_G(p)}{p-p_c}<\infty,
\end{equation}
then
\[
\Psi_{p_c}(\neg{\cal C}_t \mbox { occurs for every $t$})=1.
\]

Our goal here is to find sharper conditions of this kind in the case
where $G$ is a tree, which with the currently available
technology is virtually the only case that
can be handled with some precision. An ultimate goal might be to find some
``if and only if'' criterion in terms of 
the near-critical behaviour of $\theta_G(p)$, determining
whether there are percolating times for the critical dynamical percolation
process. However, we shall see in Remark 2.4 below 
that a general criterion of this kind is impossible, at least if one insists
that faster growth of $\theta_G(p)$ should make it easier to get percolating
times in the dynamical percolation process. 

Let $\Gamma$ be an infinite locally finite tree with vertex set $V$, edge
set $E$, and a distinguished vertex $\rho\in V$ called the root. For
$v\in V$, write $|v|$ for the distance between $v$ and $\rho$. If $v$ and 
$w$ are nearest neighbours with $|w|=|v|+1$, then $w$ is called a child 
of $v$, and $v$ is called a parent of $w$. 
The vertex set $\{v\in V:\, |v|=n\}$ is called the $n$:th level
of $\Gamma$, and is denoted $\Gamma_n$. Further, $|\Gamma_n|$ denotes
the number of vertices on the $n$:th level. If two vertices at the same level
always have the same number of children, then $\Gamma$ is called 
spherically symmetric.

Restricting to spherically symmetric trees, we do have a sharp criterion:

\medskip\noindent
{\bf Theorem 1.1:}
{\sl
Let $\Gamma$ be a spherically symmetric tree with $p_c(\Gamma)=p_c\in(0,1)$
and $\theta_\Gamma(p_c)=0$. Then 
\begin{equation} \label{no_percolating_times}
\Psi_{\Gamma,p_c}(\rho \mbox{ never percolates})=1
\end{equation}
if and only if
\begin{equation} \label{divergent_integral}
\int_{p_c}^1\frac{1}{\theta_\Gamma(p)}dp =\infty.
\end{equation}
}

\medskip\noindent
Note that (\ref{divergent_integral}) is a strictly weaker condition than
(\ref{derivative}): An example of a spherically symmetric tree $\Gamma$ 
which satisfies (\ref{divergent_integral}) but not (\ref{derivative}) can
be obtained by letting $|\Gamma_n|$ be of the order $2^n\log n$ for each $n$.
This can be shown using Lemma 2.1 below. 

Only half of Theorem 1.1 can be extended to general trees:

\medskip\noindent
{\bf Theorem 1.2:}
{\sl
Let $\Gamma$ be an infinite locally finite tree with $p_c(\Gamma)=p_c\in(0,1)$
and $\theta_\Gamma(p_c)=0$. Then {\rm (\ref{divergent_integral})} implies 
{\rm (\ref{no_percolating_times})}. 
The converse is not true, i.e.\ there exists
a choice of $\Gamma$ such that {\rm (\ref{no_percolating_times})} holds and
\begin{equation} \label{convergent_integral}
\int_{p_c}^1\frac{1}{\theta_\Gamma(p)}dp <\infty.
\end{equation}
}

\medskip\noindent
Our main motivation for this study was to 
try to shed some light on the important
case where $G$ is the square lattice ${\bf Z}^2$ (more precisely,
$V={\bf Z}^2$ and $E$ consists of pairs of vertices at Euclidean distance
$1$ from each other). It is a classical result of Kesten \cite{K80} that
$p_c=\frac{1}{2}$ and $\theta_{G}(p_c)=0$ for $G={\bf Z}^2$. Moreover, 
Kesten and Zhang 
\cite{KZ} showed that the condition (\ref{divergent_integral})
fails in this case. This means that if one could extend 
Theorem 1.1 to some class of graphs which includes ${\bf Z}^2$, then one would
be able to conclude that for critical 
dynamical percolation on ${\bf Z}^2$ there would be exceptional times
with infinite open clusters 
(this would contrast with the case of ${\bf Z}^d$ with
$d$ sufficiently large; see \cite{HPS}). One should not feel too
discouraged by Theorem 1.2 in taking up this line of research, because
the counterexample used to show that the ``only if'' direction of
Theorem 1.1 does not hold in the generality of Theorem 1.2 is highly
irregular and nonsymmetric. Theorem 1.1 thus provides some weak evidence
that ${\bf Z}^2$ might have exceptional times at criticality. We want to
stress, however, that we still do not think it is clear what the right
conjecture should be. 
A different approach to the problem of critical dynamical percolation on
${\bf Z}^2$ is discussed by Benjamini, Kalai and Schramm \cite{BKS}. 

\section{Proofs}

We first need some terminology on flows and electrical networks
on trees; see e.g.\ Lyons and Peres
\cite{LP} for a general introduction to this subject. 
An edge connecting two vertices $v,w\in V$ is denoted $\langle v,w \rangle$. 
A unit flow ${\cal F}$ on the tree $\Gamma$ is an assignment of non-negative
numbers $\{{\cal F}(e)\}_{e\in E}$ to the edges of $\Gamma$, satisfying
\begin{description}
\item{(i)} $\sum_v{\cal F}(\langle \rho, v\rangle)=1$, where the sum runs over
all children of the root $\rho$, and
\item{(ii)} for each vertex $v$ (not equal to $\rho$) with parent $w$ 
and children $v_1,\ldots, v_k$, we have 
$\sum_{v_i}{\cal F}(\langle v, v_i \rangle)={\cal F}(\langle w,v \rangle)$. 
\end{description}
This should be thought of as having ``flow in equal to flow out'' in all
vertices except $\rho$ whose net flow out is $1$. An electrical network
${\cal C}$ on $\Gamma$ is simply an assignment of arbitrary positive numbers
$\{{\cal C}(e)\}_{e\in E}$, called
conductances, to the edges of $\Gamma$. The energy $W$ 
of a unit flow ${\cal F}$ in the network ${\cal C}$ is defined as
\[
W({\cal F}, {\cal C}) =
\sum_{e\in E} \frac{({\cal F}(e))^2}{{\cal C}(e)} ,
\]
and the effective conductance of the network ${\cal C}$ is defined
as $[\inf\{W({\cal F}, {\cal C})\}]^{-1}$ where the infimum runs over all
unit flows ${\cal F}$
on $\Gamma$. This infimum is in fact a minimum, which moreover
is attained at a unique unit flow whenever it is nonzero. 
Let $C(\Gamma,p)$ denote the
effective conductance of the electrical network
$N_{\Gamma,p}$ obtained by assigning each edge
between $\Gamma_{n-1}$ and $\Gamma_n$ conductance
$p^n$. Similarly, let $C^\ast(\Gamma,p)$ 
be the effective conductance the network $N^\ast_{\Gamma,p}$ in which
edges between levels $n-1$ and $n$ are assigned conductance $np^n$.  

The keys to the proof of 
Theorem 1.1 are the following two results from the literature. 

\medskip\noindent
{\bf Lemma 2.1 (Lyons \cite{L}):} 
{\sl For each tree $\Gamma$ and each $p\in(0,1)$ we have
\begin{equation} \label{eq:Lyons}
\frac{C(\Gamma,p)}{1+C(\Gamma,p)}\leq \theta_\Gamma(p)\leq 2
\frac{C(\Gamma,p)}{1+C(\Gamma,p)}.
\end{equation}
}

\medskip\noindent
(See Marchal \cite{M} for a recent sharpening of the upper bound in 
(\ref{eq:Lyons}).)

\medskip\noindent
{\bf Lemma 2.2 (H\"aggstr\"om, Peres and Steif \cite{HPS}):} 
{\sl For any tree $\Gamma$ with critical value $p_c$ we have that
{\rm (\ref{no_percolating_times})} holds if and only if
\[
C^\ast(\Gamma, p_c)=0.
\]
}

\medskip\noindent
{\bf Proof of Theorem 1.1:}
Since $\Gamma$ is spherically symmetric, we have
\begin{equation} \label{ss_conductance}
C(\Gamma,\rho)=\left( \sum_{k=1}^\infty \frac{p^{-k}}{|\Gamma_k|}\right)^{-1}
\end{equation}
and
\begin{equation} \label{ss_star_conductance}
C^\ast(\Gamma,\rho)=
\left( \sum_{k=1}^\infty \frac{p^{-k}}{k|\Gamma_k|}\right)^{-1}.
\end{equation}
By Lemma 2.1, condition (\ref{divergent_integral}) can be rewritten as
\[
\int_{p_c}^1\frac{1+C(\Gamma,p)}{C(\Gamma,p)}\, dp=\infty
\]
which obviously is equivalent to
\begin{equation} \label{int_resistance}
\int_{p_c}^1\frac{1}{C(\Gamma,p)}\, dp=\infty.
\end{equation}
By (\ref{ss_conductance}), this is the same as
\[
\int_{p_c}^1 \sum_{k=1}^\infty \frac{p^{-k}}{|\Gamma_k|}\, dp=\infty 
\]
which, in turn, is the same as
\begin{equation} \label{sum_of_integrals}
\sum_{k=1}^\infty\frac{\int_{p_c}^1 p^{-k}\, dp}{|\Gamma_k|}=\infty
\end{equation}
by Fubini's Theorem. The left hand side equals
\[
\frac{\ln(\frac{1}{p_c})}{|\Gamma_1|}+\sum_{k=2}^\infty
\frac{p_c^{-(k-1)}-1}{(k-1)|\Gamma_k|}
\]
so that (\ref{sum_of_integrals}) holds if and only if
\[
\sum_{k=1}^\infty\frac{p_c^{-k}}{k|\Gamma_k|}=\infty.
\]
Now (\ref{ss_star_conductance}) and Lemma 2.2 complete the proof. $\Cox$

\medskip\noindent
For the proof of Theorem 1.2, it is convenient to isolate the following 
lemma.  

\medskip\noindent
{\bf Lemma 2.3:}
{\sl 
For any $p,q\in(0,1)$, there exists a spherically symmetric tree $\Gamma$ with
$p_c(\Gamma)=p$ and $\theta_\Gamma(p)\geq q$. 
}

\medskip\noindent
{\bf Proof:} Fix $p$ and $q$. Let $\Gamma^\prime$ be a spherically symmetric
tree with $|\Gamma^\prime_k|$ being bounded above and below by constants times
$k^2p^{-k}$. Such a tree can be defined inductively by letting each vertex
on level $k-1$ have exactly
\[
\min\{i\in\{1,2,\ldots\}:\, i|\Gamma^\prime_{k-1}|\geq k^2p^{-k}\}
\]
children. It is an easy application of Lemma 2.1 
to show that $p_c(\Gamma^\prime)=p$ and
$\theta_{\Gamma^\prime}(p)>0$. If $\theta_{\Gamma^\prime}(p)\geq q$ we are
done by taking $\Gamma=\Gamma^\prime$. Otherwise we set 
\[
j=\min\{i\in\{1,2,\ldots\}:\, (1-\theta_{\Gamma^\prime}(p))^j\leq 1-q\}.
\]
and let $\Gamma$ consist
of $j$ copies of $\Gamma^\prime$ sharing the same root $\rho$. 
$\Cox$

\medskip\noindent
{\bf Proof of Theorem 1.2:}
We first show that (\ref{divergent_integral}) implies 
(\ref{no_percolating_times}) by proving the equivalent statement that
\begin{equation} \label{perc_at_criticality}
\Psi_{p_c}(\exists t\mbox{ at which $\rho$ percolates})>0
\end{equation}
implies (\ref{convergent_integral}). 
Suppose that (\ref{perc_at_criticality}) holds. We then have 
$C^\ast(\Gamma,p_c)>0$ by Lemma~2.2. Let ${\cal F}^\ast$ 
denote the minimal energy 
unit flow from $\rho$ to infinity in the electrical network 
$N^\ast_{\Gamma,p_c}$. Write $W^\ast$ for the energy of ${\cal F}^\ast$ in
$N^\ast_{\Gamma,p_c}$, so that 
\begin{equation} \label{finite_energy}
W^\ast=\frac{1}{C^\ast(\Gamma,p_c)}<\infty.
\end{equation}
Decompose $W^\ast$ as
\begin{equation} \label{decomposition}
W^\ast=\sum_{k=1}^\infty W^\ast(k)
\end{equation}
where $W^\ast(k)$ is the contribution to $W^\ast$ coming from edges connecting
levels $k-1$ and $k$ in $\Gamma$. Also write $W_p$ for the energy of 
${\cal F}^\ast$
in the network $N_{\Gamma,p}$, and decompose $W_p$ into
\[
W_p=\sum_{k=1}^\infty W_p(k)
\]
analogously to (\ref{decomposition}). We have
\[
\frac{1}{C(\Gamma,p)}\leq W_p ,
\]
and furthermore
\[
W_p(k)=k\left( \frac{p_c}{p}\right)^k W^\ast(k).
\]
Hence
\begin{eqnarray*}
\int_{p_c}^1\frac{1}{C(\Gamma,p)}\, dp & \leq & 
\int_{p_c}^1 \sum_{k=1}^\infty W_p(k)\, dp \\
& = & \sum_{k=1}^\infty \int_{p_c}^1 W_p(k)\, dp \\
& = & \sum_{k=1}^\infty kW^\ast(k) 
\int_{p_c}^1 \left( \frac{p_c}{p} \right)^k \, dp \\
& = & p_c\sum_{k=1}^\infty \frac{k(1-p_c^{k-1})}{k-1}W^\ast(k) \\
& < & \infty
\end{eqnarray*}
where the last step is due to (\ref{decomposition}) and 
(\ref{finite_energy}). The desired
conclusion (\ref{convergent_integral}) follows using Lemma 2.1. 

It remains to give an example of a tree $\Gamma$ for which 
(\ref{no_percolating_times}) and (\ref{convergent_integral}) both hold. 
For each $p\in(1/2,1)$, let
$\Gamma^{(p)}$ denote some tree with critical value $p_c(\Gamma^{(p)})=p$
and $\theta_{\Gamma^{(p)}}(p)\geq 1/2$; such trees exist by Lemma 2.3. 
The construction of $\Gamma$ is as follows. Start with a single
infinite branch $(\rho, v_1, v_2, \ldots)$ and attatch to each vertex 
$v_i$ a copy of $\Gamma^{(p_i)}$, where the sequence $\{p_i\}_{i\geq 1}$
is decreasing with $p_i>1/2$ for each $i$ and 
$\lim_{i\rightarrow\infty} p_i =1/2$. We get $p_c(\Gamma)=1/2$. Furthermore,  
$\theta_\Gamma(p_i)\geq 2^{-(i+1)}$, so that we can make $\theta_\Gamma(p)$ 
take off from
criticality arbitrarily fast by letting the $p_i$'s tend to 1/2 sufficiently
fast. In particular, we can make (\ref{convergent_integral}) 
hold by e.g.\ taking
$p_i=1/2 + 3^{-i}$. 
On the other hand, note that percolation on $\Gamma^{(p_i)}$ is subcritical
at $p=1/2$ for each $i$, whence
\[
\Psi_{\Gamma^{(p_i)},1/2}[ \forall t \neg{\cal C}_t]=1.
\]
Since the existence of an infinite open cluster in $\Gamma$ implies the
existence of an infinite cluster either in the branch 
$(\rho, v_1, v_2,\ldots)$ or in one of the subtrees $\Gamma^{(p_i)}$, we get
\[
\Psi_{\Gamma,1/2}[ \forall t \neg{\cal C}_t]=1,
\]
so the proof is complete. 
$\Cox$

\medskip\noindent
{\bf Remark 2.4:} The construction at the end of the proof of Theorem 1.2 shows
that we can find two trees $\Gamma$ and $\Gamma^\prime$ with the properties
that
\begin{description}
\item{(i)}
$p_c(\Gamma)=p_c(\Gamma^\prime)=1/2$,
\item{(ii)}
$\theta_\Gamma(1/2)=\theta_{\Gamma^\prime}(1/2)=0$,
\item{(iii)}
the critical dynamical percolation process has times with infinite clusters
for $\Gamma$ but not for $\Gamma^\prime$, and
\item{(iv)}
there exists a $p^\ast>1/2$ such that
$\theta_\Gamma(p)<\theta_{\Gamma^\prime}(p)$ for all $p\in (1/2,p^\ast)$. 
\end{description}
Indeed, take $\Gamma$ to be spherically symmetric with $|\Gamma_n|$ 
of the order $n2^n$, and then build up $\Gamma^\prime$ as in the above
construction sending $p_i$ to $1/2$ fast enough so that (iv) holds. 
In words, (i)--(iv) tell us that
$\Gamma^\prime$ is closer to percolating at criticality than
$\Gamma$ in the sense of how fast $\theta(p)$ takes off from criticality,
whereas $\Gamma$ is closer to percolating in the sense of dynamical 
percolation.

\end{document}